\documentclass[11pt]{pjm}

\usepackage{amssymb, amsmath, latexsym, amscd}
\usepackage{epsfig}
\usepackage{psfrag}
\newtheorem{theorem}{Theorem}

\newtheorem{definition}{Definition}

\newtheorem{corollary}{Corollary}

\newtheorem{proposition}{Proposition}

\newcommand{\G}{\mathcal G}
\newcommand{\g}{\mathcal G}

\renewcommand{\O}{\mathcal O}

\title{Fusion and fission in graph complexes}

\author[J. Conant]{James Conant}
\address{Dept. of Mathematics\\
         Cornell University\\
         Ithaca, NY 14853-4201}
\email{jconant@math.cornell.edu}
\urladdr{http://www.math.cornell.edu/\~{}jconant/pageone.html}
  
\thanks{Partially supported by NSF VIGRE grant DMS-9983660}

\theoremstyle{definition}
  
\date{Received date / Revised version date}
 
\begin{document}
 
\begin{abstract}
 We analyze a functor from cyclic operads to chain complexes 
first considered by Getzler and Kapranov and 
also by Markl. This functor is a generalization of the graph homology
considered by Kontsevich, which was defined for the three operads ${Comm},
Assoc$, and ${Lie}$. More specifically we show that
these chain complexes have a rich algebraic structure 
in the form of families of operations defined by 
\emph{fusion} and
\emph{fission}. These operations fit together to form uncountably
many $Lie_\infty$ and $co$-$Lie_\infty$ structures. 
  In particular, the chain complexes have a
 bracket and cobracket which are compatible in the Lie
bialgebra sense on a certain natural subcomplex. 
\end{abstract}
 
\maketitle
 
\section{Introduction}
More than a decade ago Maxim Kontsevich \cite{kontsevich} considered
graph homology as a tool for studying and computing the homology of many
seemingly disparate objects. 
One version of the graph complex computes, via work of R.C. Penner
\cite{penner}, the homology of the moduli space (or equivalently mapping
class group) of surfaces. Another version computes, via work of M. Culler
and K. Vogtmann \cite{culler-vogtmann},  the homology of the group of
outer automorphisms of the free group. There is also a version which is
related to finite type invariants of three-manifolds. On the other hand,
these three graph complexes compute the homology of three infinite
dimensional Lie algebras, leading to quite unexpected isomorphisms.
Kontsevich's graph complexes were generalized to the case of
\emph{modular operads} by Getzler and Kapranov\cite{gk2}, and were 
considered for the special case
 of cyclic operads
by Martin Markl \cite{markl}.

In \cite{cv} Karen Vogtmann and I showed that the commutative
graph complex carries the structure of both a Lie algebra and a Lie
coalgebra. These are compatible as a bialgebra on a certain natural
subcomplex.
In this paper I will generalize these two operations to the case of any
cyclic operad, and show that they are each first in a series of
higher order operations which fit together nicely and vanish on homology.

Let the graph complex corresponding to a cyclic operad $\O$ be denoted by
$\g^\O$.
I will define a sequence of ``higher order brackets"
$$\phi_n\colon S^n\g^\O\to \g.$$
The map $\phi_n$ is defined by fusing together $n$ graphs along a
$2n$-gon in all possible ways (Figure~\ref{ngon}).
Extending each $\phi_n$ as a coderivation to $S\g^\O$,
 these maps are all compatible with each other in a very
strong sense (Theorem~\ref{lemma1}). For any subset $I\subset \mathbb N$,
define
$\phi_I =\sum_{i\in I}\phi_i$. Theorem~\ref{lemma1} implies that
$\phi_I^2=0$. This is precisely the definition of
 a $Lie_\infty$ (strong homotopy Lie) structure.
 In this way we get
uncountably many $Lie_\infty$ structures.

Let $P\g^\O$ denote the subcomplex of the graph complex spanned by 
connected graphs. I will define a sequence of ``higher
order cobrackets"
$$\theta_n: P\g^\O\to  S^nP\g^\O.$$
The map $\theta_n$ is defined by fissioning
a graph into $n$ graphs along a $2n$-gon (Figure~\ref{ngon2}).
The $\theta_n$ maps, extended to $S\G^\O$ as derivations, are compatible
in a strong sense also (Theorem 2).  For any $I\subset\mathbb N$, 
$\theta_I$ is defined as above, and 
Theorem 2 implies each $\phi_I$ is a $co$-$Lie_\infty$
structure.

Trouble arises, as
was foreshadowed in \cite{cv} in the compatibility between brackets
and cobrackets. In \cite{cv} we were able to avoid difficulty by
restricting to connected graphs without separating edges,
and indeed in this context $\theta_2,\phi_2$ are compatible in a Lie
bialgebra sense (Theorem~\ref{compat}).  
But there
appears to be no similar way out for the higher order operations. 
The higher order brackets and cobrackets simply
fit
together in a more complicated way than one would guess, even on graphs
without separating edges.

All of the operations
are highly nontrivial on chains, and are compatible with the boundary
operator. Indeed they vanish canonically on the level of homology. Thus
these operations can be thought of as ``generalized Schouten brackets,"
since in the case of Lie algebras, the Schouten bracket is an operation
on the Chevalley-Eilenberg complex which vanishes canonically upon
application of the homology functor.

 Moira Chas and Dennis Sullivan \cite{cs} define similar structures
on
\emph{string homology}, the homology of a free loop space. They define an
uncountable family of
$Lie_\infty$ structures, indexed by sets of positive integers, on string
homology which obey the same compatibility relations as the ones found
here (Theorem 1). They also find a Lie bialgebra structure
\cite{chas},\cite{cs2}.  Drawing the analogy further, one is led to
speculate that string homology has an uncountable infinity of
$co$-$Lie_\infty$ structures. It would be interesting to know whether such
co structures, if they exist, are compatible in a nice way with the Lie
structures, or if they mirror more complicated graph interactions.  

{\bf Acknowledgements:} It is a pleasure to thank Karen Vogtmann for many
discussions. I also wish to thank Swapneel Mahajan for his perceptive
input. Credit also goes to the anonymous referee who noticed an error in
the original manuscript and suggested many expositional improvements.
\newpage

\section{Cyclic operads and graphs}
We begin by briefly reviewing the salient features of a cyclic operad,
and proceed to give
Markl's construction of graph complexes. A good introduction
to these objects can be found in the recent book by Markl, Shnider and
Stasheff \cite{mss}. 

Kontsevich's three graph complexes are associated to the commutative,
associative and Lie operads. Each of these operads
${\mathcal O}=\oplus \mathcal O(n)$ has a description as a vector space
spanned by different flavors of rooted trees with labelled leaves.

\begin{figure}
\begin{center}
  \epsfig{file=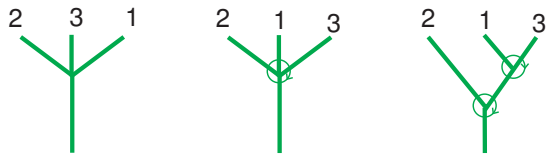,height=2cm}
  \caption{Elements of $Comm[3],Assoc[3],$ and $Lie[3]$ respectively.}
  \label{operads}
\end{center}
\end{figure}

\begin{figure}
\begin{center}
\epsfig{file=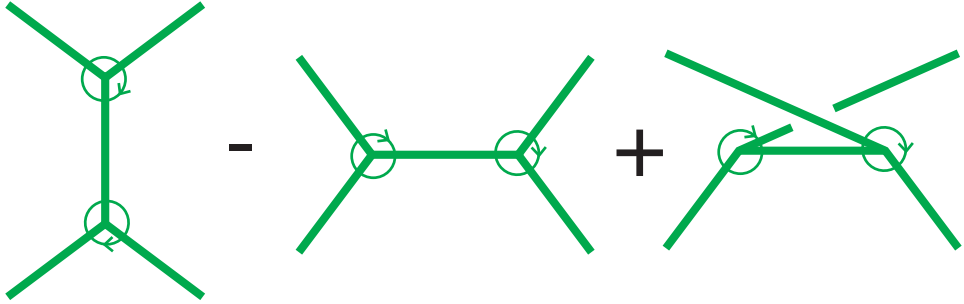, height=4cm}
\caption{The IHX relation. Each term represents a piece of a graph which
is identical outside of the pictured spot.\label{ihxagain}}
\end{center}
\end{figure}

\begin{itemize}
\item The $n$th degree part of the commutative operad $Comm(n)$ has a
basis consisting of rooted trees which have one internal vertex and
$n$ labelled leaves. Hence $Comm[n]$ is $1$ dimensional! 
The composition law $$Comm[m]\otimes Comm[n_1]\otimes\ldots\otimes
Comm[n_m]\to Comm[n_1+\ldots +n_m]$$ is defined on $c\otimes
c_1\otimes\ldots\otimes c_m$ by grafting the root of each $c_i$ onto the
leaf of $c$ labelled by $i$ for each $i$, and suitably relabelling the
leaves. The composition is completed by contracting all edges not
adjacent to a root or a leaf.

\item The $n$th degree part of the associative operad $Assoc(n)$ has a
basis consisting of rooted trees with one internal vertex which have a
specified cyclic ordering of the edges incident to the vertex, and which
have $n$ labelled leaves. Composition is again by grafting and
contracting created edges, with the proviso that the cyclic ordering is
respected.

\item The $n$th degree part of the Lie operad $Lie(n)$ is actually easiest
to describe as a quotient space $\overline{Lie(n)}/AS + IHX$.
$\overline{Lie(n)}$
 has a basis given by rooted trivalent trees with $n$ labelled
leaves, where each vertex has a specified cyclic order of adjacent edges.
The AS subspace is spanned by sums $T_1+T_2$, where $T_{1,2}$ are
identical except for a cyclic ordering on some vertex. Modding out by AS
says that Lie algebras are anti-symmetric.
 The IHX subspace is spanned by sums $T_1-T_2+T_3$, where $T_{1,2,3}$ are
identical trees except at one spot where they are as in Figure
\ref{ihxagain}. On the level of Lie algebras this is the Jacobi relation.
Composition is via grafting, but without the contraction step. 
\end{itemize}

Notice that in each of these cases the action of the symmetric group
$S_n$ which permutes the labels of the leaves can be extended to an
action of $S_{n+1}$. This is by thinking of the root as another labelled
leaf, say labelled by $0$. (One must check in the Lie case that the IHX
subspace is preserved by this action.) Operads where this extension is
possible are called \emph{cyclic} \cite{gk}, provided that
the extension satisfies appropriate axioms. 
Other examples of cyclic
operads are the endomorphism operad and the Poisson operad.

As a general philosophy, one can think of cyclic operads as consisting of
unrooted trees, with composition given by some version of grafting.
 The idea is to plug these in to the nodes of a
graph to obtain different types of decorations on a graph. Plugging in
a basis element from $Comm(n)$ at each vertex of valence $n$, one simply
gets an undecorated graph. Plugging in an element of $Assoc(n)$ one gets a
cyclic order at the vertex. This is often called a ribbon graph.
Plugging in an element of $Lie(n)$ gives a relatively strange object.
By definition it is obtained from some unrooted ribbon trivalent
trees by joining the leaves together with edges. 
See Figure \ref{Lieexample}.
Thus one may think of it as a
trivalent graph with a special distinguished subset where IHX and AS
relations may take place.
 It
is reminiscent of the diagram algebras that appear in the theory of
Vassiliev invariants of low dimensional objects, which consist of
(uni)trivalent graphs, but where the AS and IHX relations are not
restricted to a distinguished subset.     

In general, let $H(v)$ be the set of half-edges incident to a vertex.
Let $L$ be a labelling of the elements of $H(v)$ by $0,\ldots,n$, where
$n+1$ is the valence of $v$.  Now define 
$$\mathcal O ((H(v))) = \left(\oplus_L \mathcal O(n)\right)_{S_{n+1}}$$
which is the set of coinvariants under the action of $S_{n+1}$, which
acts as follows. If $o\in \mathcal O(n)$, let
$o_L$ denote putting $o$ in the $L$th summand of the direct sum.
If $\sigma\in S_{n+1}$ then $\sigma \cdot o_L =
(\sigma\cdot o)_{\sigma\cdot L}$. When $\O$ is an operad of trees,
$\O ((H(v)))$ is isomorphic to the space of identifications of the leaves
and root of elements of $\O[n]$ with the half-edges incident to $v$.

Now we define an $\mathcal O$-labelling of a graph to be a choice of
element $o_v\in\mathcal O ((H(v)))$ for each vertex, $v$, of the graph.
Graphically, we put a circle at each vertex to represent the operad
element.

In addition we would like a notion of ``orientation" of a graph, which
will make it possible to define a boundary operator. This is analogous to
the need for an orientation of the simplices of a simplicial complex in
order to do the same. There are many equivalent notions,
perhaps the most intuitive is the following. 

\begin{definition}
An orientation of a graph is an ordering of the
vertices and a choice of direction for each
edge, modulo the even action of $S_V\times \mathbb Z_2^E$. 
Here $V$ and $E$ are the number or vertices and edges of the graph,
respectively. An element of
$S_V\times \mathbb Z_2^E$ is called even if it is a product of an even
number of elements each of which is either a transposition in $S_v$ or an
element of the form
$(0,\ldots,1,\ldots,0)\in \mathbb Z_2^E$.
\end{definition}

Notice that any graph has exactly two orientations. Let $-$ indicate the
map switching orientations.
\vskip10pt
{\bf Remark:} Lie graphs actually have a much simpler description, because
the orientations of the graph and vertices cancel out to a large degree.
Namely, one can think of a Lie graph as a trivalent graph with a
distinguished subforest, whose edges are ordered modulo even permutations.
The IHX relation in the Lie operad becomes the condition that 
the three terms in an IHX relation of the trivalent graph sum to
zero provided the edge involved is in the forest. This will be explained
carefully in
\cite{cv2}.
 
\subsection{Chain complexes}
Now for any cyclic operad $\mathcal O$ we are ready to define $\mathcal
O$-graph complexes.

Define $\mathcal G^{\mathcal O}_v$ to be the span of $\mathcal
O$-labelled oriented graphs with vertices all of valence $\geq 3$,
modulo the relation $(G,or)=-(G,-or)$ and also modulo multilinearity of
the $\O$-labels. More precisely, we set $$\mathcal G^{\mathcal
O}=\left({\bigoplus_{(G,or)}
\bigotimes_{v\in V(G)}\mathcal O
((H(v)))}\right){\text{\Huge /}}\begin{array}{c}{ }\\
\{(G,or)=-(G,-or)\}\end{array} ,$$ where the
direct sum is over oriented graphs with vertices of valence $\geq 3$, and
where
$V(G)$ is the set of vertices of $G$.
Define $\mathcal G^{\mathcal O}_v$ to 
be the part of $\g^\O$ with $v$ vertices.

\begin{figure}
\begin{center}
\epsfig{file=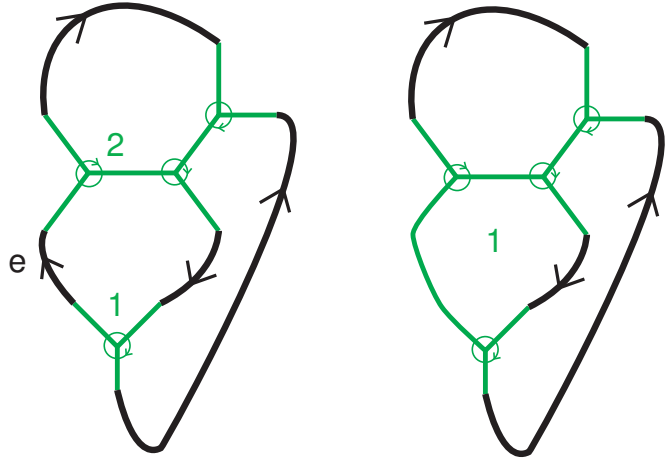, height=7cm}
\caption{The Lie graphs $G$ and $G_e$.\label{Lieexample}} 
\end{center}
\end{figure}

For each edge, $e$, in a graph $(G,or)$ we define contraction along that
edge $(G,or)_e$ to be the graph where the two operad elements at each
endpoint of $e$ are composed along $e$. The induced orientation can be
fixed by assuming that the endpoints of $e$ are labelled $1$ and $2$ and
the edge direction is from $1$ to $2$. The new vertex, which results
from composing the two operad elements, is labelled $1$, and all other
indices are reduced by $1$. If $e$ is a loop, then define $(G,or)_e=0$.
In the commutative case, $(G,or)_e$ is defined by simply contracting the
edge of the (undecorated) graph. In the associative case the cyclic
orders at both endpoints of an edge are joined together to give a cyclic
order at the vertex resulting from the edge collapse.
For an example in the Lie case, see Figure \ref{Lieexample}.

Define $\partial_G:\G^{\mathcal O}_v\to \G^{\mathcal O}_{v-1}$ by
$\partial_G(G,or)=\sum_{e\in E(G)} (G,or)_e$, where $E(G)$ is the set of
edges of $G$.
\vskip6pt
{\bf Remark:} In the simpler version of the Lie case,
 the boundary operator adds
an edge to the forest in all possible ways, with the edge's number coming
directly after the edge numbers in the original forest.

\vskip10pt
$\G^{\mathcal O}$ is a graded commutative algebra under disjoint union.
It is also a graded commutative coalgebra, defining the coproduct such
that connected graphs are primitive, and extending multiplicatively.
Thus we may write $P\G^{\mathcal O}$ for the subspace generated by
connected graphs.
Let $P^{(n)}\G^{\mathcal O}$ be the subspace generated by
connected graphs with $b_1 = n$. 

Even though this paper is concerned with chain complexes and not their
homology per se, it is still useful to record the following facts.

Let $Out(F_n)$ denote the group of outer automorphisms of the free group
$F_n$, and let $\mathcal M^\prime_{g,m}$ denote the moduli space of a
surface of genus $g$ with $m$ unlabelled punctures.

Then 
\begin{align*}
H_k(P^{(n)}\G^{\mathcal Assoc}) &= 
\bigoplus_{m\geq 1,g:2g+m-1=n}H^{4g-4+2m-k}(\mathcal M^\prime_{g,m};
\mathbb Q)\\
 H_k(P^{(n)}\G^{\mathcal Lie}) &= H^{2n-2-k}(Out(F_n);\mathbb Q)
\end{align*} 
In addition, part of commutative graph cohomology plays a role in the
theory of finite type invariants of homology $3$-spheres. More precisely,
 we have that 
$$\bigoplus_{n\geq 2} H^{2n-2}(P^{(n)}\G^{\mathcal Comm}) \cong P\mathcal
A(\emptyset)$$
where $P\mathcal A(\emptyset)$ is the diagram algebra where the
logarithm of the Aarhus version of the LMO invariant takes
values \cite{bngrt}. 

The first two statements above are due to Kontsevich, being implicit in
his paper \cite{kontsevich}. A more detailed explanation of these two
facts and their proofs will appear in \cite{cv2}. The last statement, the
relation to finite type invariants, is essentially content-free, being a
trivial isomorphism, at least modulo equivalences of various notions of
orientation. 

\subsection{Cohomology}
In at least two interesting cases, it is possible to define graph
cohomology. The coboundary operator $\delta$ is the sum of
inserting an edge in all possible ways. In the commutative and
associative cases this makes perfect sense. Unfortunately, in the Lie
case an insertion, which is essentially the deletion of an edge from
the forest, does not preserve the IHX subspace and is not well-defined.
 In the cases where $\delta$ can be defined the boundary and
coboundary are adjoint with respect to the inner product $<G,H> =
|Aut(G)|\delta_{GH}$. This can be seen by applying the argument of
\cite{cv}, Proposition 12 mutatis mutandis.

\section{Fusion}

\begin{figure}
\begin{center}
\epsfig{file=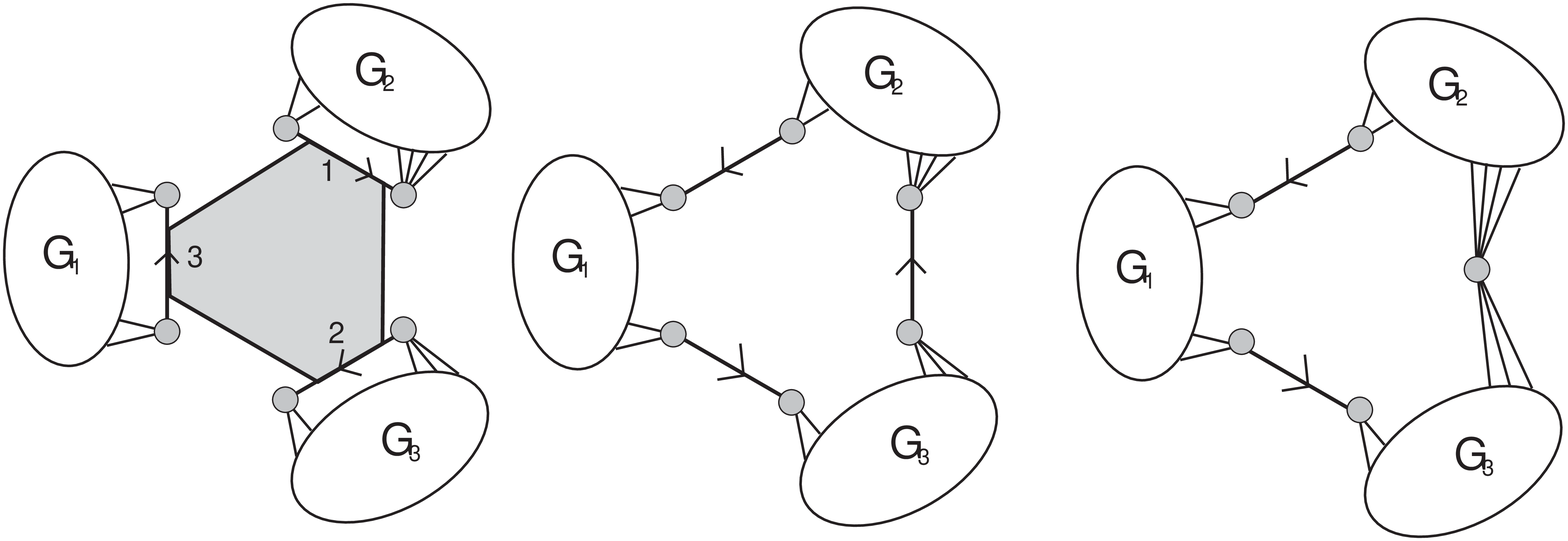,height=4cm}
\caption{One term in $\phi_3(G_1,G_2,G_3)$.\label{ngon}}
\end{center}
\end{figure}

We start with an oriented labelled $2n$-gon. Label every other edge on its
perimeter consecutively by the numbers $1\ldots n$, consistent with the
orientation. Now fix $n$ directed edges $e_1,\ldots, e_n$ of a graph $G$.
Define $G~<~e_1,\ldots, e_n~>$ to be the graph formed in the
following way. First, for each $i$, glue the edge marked $i$ of the
$2n$-gon to the edge $e_i$ of the graph. Second, delete these edges along
which the $2n$-gon was attached leaving $n$ new edges. This is
illustrated in Figure \ref{ngon}. The graph $G<e_1,\ldots, e_n>$ has an
induced orientation which can be easily described. Fix a labelling of the
graph such that the directions of the edges $e_1,\ldots,e_n$ are both
consistent with the graph's orientation and with the directions which
correspond to the gluing. The $n$ new edges have orientations induced by
the $n$-gon. Switch all of these, as is usual with a cobordism.
Now, for any $n\in \mathbb N$ we define an operation 
$$\phi_n : S^n \g^\O \to \g^\O$$
by $\phi_n(G_1\odot\cdots\odot G_n) = \sum (G_1\cdot G_2\cdots\cdot
G_n)<e_1,\ldots,e_n>_e$, where 
\begin{itemize}
\item The sum is over all $n$-tuples of directed
edges
$(e_1,e_2,\ldots,e_n)$ all of which Lie in separate $G_i$. 
\item The notation ``$\odot$" denotes ``graded symmetric tensor product."
\item The edge $e$
which is contracted is the edge coming from the boundary of the $2n$-gon
between ``1'' and ``2.'' 
\end{itemize}
Thus $\phi_n$ is a type of fusion operation which
takes
$n$ graphs and fuses them together along a $2n$-gon.

Extend $\phi_n$ to $S\g^\O$ as a coderivation. That is
$$\phi_n(G_1\odot\ldots\odot G_k) = \sum_{I\cup J} {\epsilon
(I,J)} \phi_n ( G_I ) \odot
G_J,$$
where $I,J$ is an ordered partition of $1,\ldots , k$, with $|I|=n$, and
$\epsilon(I,J)$ is the sign defined by the equation $G_1\odot\ldots\odot
G_k = \epsilon(I,J)G_I\odot G_J$. Notice that 
 $\phi_1$ by definition glues on a bigon to an edge, which doesn't change
the edge, and then contracts it. That is, $\phi_1 = \partial_G$. Notice
that it doesn't matter whether we extend $\partial_G$ to $S\G^\O$ as a
derivation or a coderivation, since they are equivalent in this case!

\begin{theorem} \label{lemma1}
The following equations hold.
\begin{itemize}
\item[a)]$\forall i \hskip5pt \phi^2_i = 0$
\item[b)]$\forall i\neq j \hskip5pt \phi_i\phi_j+\phi_j\phi_i=0$. 
\end{itemize}
\end{theorem} 

\begin{corollary}
Let $I$ be a subset of $\mathbb N$, finite or infinite. Let
$\phi_I = \sum_{i\in I}\phi_i$. Then $\phi_I^2=0$. 
\end{corollary}

\vskip5pt

\noindent\emph{[Proof of Theorem~\ref{lemma1}]}

First we show $\phi_n^2|_{S^k\g^\O}=0$.
We only need consider the case when $k= 2n-1$, which implies the higher
degree cases. 
\begin{align*}
\phi_n^2(G_1\odot\cdots \odot G_{2n-1}) &= \sum_{I\cup J=[2n-1]}
\phi_n(\phi_n(G_I)\odot
G_J)
\end{align*}
Thus we are attaching a disk to $G_i$ where $i\in I$ along its $n$
subarcs. We then attach a disk to the result together with the other
$n-1$ graphs. If the second disk attaches to an edge not involved in the
first disk, then this gives the same unoriented result as attaching the
disks in the other order. Keeping track of the orientation, we see that
the two orders of attaching the disk cancel.
 The other possibility is that the second disk
attaches to the first. This can be thought of as attaching a $4n-2$-gon
 to the $2n-1$ graphs, with a separating arc along the
$4n-2$-gon, and two ordered edges marked for collapse.
We can simplify the combinatorics somewhat by shrinking the complement of
the $2n-1$ attaching regions for the disk, to get a $2n-1$-gon with an
arc joining two vertices and two vertices marked for collapse. The sorts
of configurations that arise are exactly recorded by the concept of
\emph{admissible} defined below.
 The lemma now
follows from the following analysis of
$2n-1$-gons.

Define $Conf(2n-1,n)$ be the set of admissible configurations of a
$2n-1$-gon. An admissible configuration consists of an embedded arc on
the $2n-1$-gon between two of the vertices, thereby partitioning the
$2n-1$ vertices into two sets of $n-1$ and $n-2$ respectively, on each
side of the arc. There are also two vertices labelled by $1$ and $2$, the
$1$ must be in the set of $n-1$ and the $2$ must be among the $n-2$ or it
could be one of the endpoints of the arc. We claim that the subset of
$Conf(2n-1,n)$ where two specific vertices are marked $1$ and $2$ is
bijective with the subset where these vertices are marked $2$ and $1$,
respectively. This follows from the fact that there is a unique
automorphism exchanging any two vertices of a $2n-1$-gon. This induces a
bijection between the two types of configurations. Keeping track of
orientations, we see that the terms of 
 corresponding to elements of $Conf(2n-1,n)$ cancel in pairs.

The fact that $\phi_i, \phi_j$ anti-commute follows from the following
similar facts about configurations of $i+j-1$-gons, $Conf(i+j-1,i,j)$.
The arc in this case will partition the vertices into a set of $i-1$, and
a set of $j-2$, where the $1$ vertex must Lie in the $i-1$ and the $2$
elsewhere. We claim there is a bijective correspondence between 
subsets of $Conf(i+j-1,i,j)$ where two fixed vertices are labelled 1 and 2
and the subsets of $Conf(i+j-1,j,i)$ where these vertices are labelled 2
and 1. To see this, fix an automorphism of the $i+j-1$-gon, exchanging
the two given vertices. This will carry one set of configurations onto
the other.
$\hfill\Box$

\begin{proposition}\label{ich}
$\phi_n$ is canonically zero at the level of homology.
\end{proposition}
\emph{[Proof]}

The fact that $\phi_n$ is even compatible with homology is the fact
$$\partial_G\circ \phi_n +\phi_n\circ \partial_G = 0$$
where $\partial_G$ is extended to $S\G^\O$ as a derivation. This follows
since $\partial_G=\phi_1$. 

It remains to show that it vanishes canonically. Consider the map
$$\mu_n\colon S^n \g^\O \to \g^\O$$ which is defined by gluing in a
$2n$-gon in all possible ways, but without contracting an edge. Then a
straightforward argument shows that $\phi_n=\partial_G\mu_n -
\mu_n\partial_G$. Thus if the input to $\phi_n$ consists of $n$ cycles,
the $\mu_n\partial_G$ term in this equation vanishes, and what is left
expresses $\phi_n$ as a boundary.$\hfill \Box$

\section{Fission} 
In this section, for simplicity, we restrict ourselves to connected
graphs, although much of it can be generalized to the nonconnected case.
In particular, when edge insertions make sense, one can dualize and
prove Theorem~\ref{whoknows} analogously to Proposition 11 of
\cite{cv}. 

Note that $\G^{\mathcal O}\cong S(P\G^{\mathcal
O})$. Denote this isomorphism by $S$. Let $$\pi_i\colon S(P\G^{\mathcal
O})\to S^i(P\G^{\mathcal O})$$ be the natural projection.
Define the map
$$\partial_i:\G_v^{\mathcal O}\to\G_{v-1}^{\mathcal O}$$
by summing over all ways of attaching a $2i$-gon to the edges of an
$\mathcal O$-graph, and then contracting the edge between $1$ and $2$. The
behavior of this operator (which does not have square zero) is
complicated, but it becomes better behaved if we look at the part which
disconnects the graph the most.

\begin{figure}
\begin{center}
   \epsfig{file=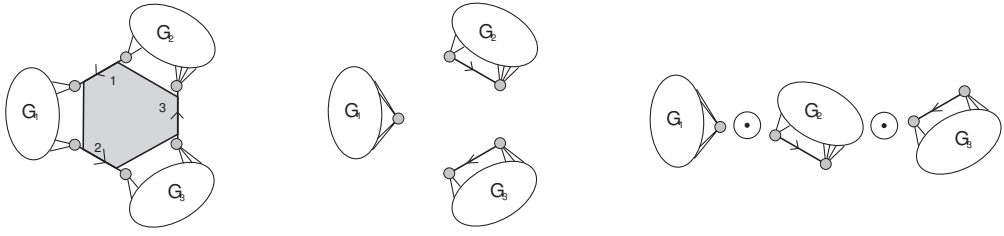,height=3cm}
   \caption{A term in $\theta_3(G)$. The middle picture
            represents a term in $\partial_3(G)$, and the final picture
            is a result of applying $S$.}\label{ngon2}
\end{center}
\end{figure}

\begin{definition}
The map $$\theta_i:P\G^{\mathcal O}\to S^i(P\G^{\mathcal
O})$$ is defined as the composition $\frac{1}{2}\pi_i\circ S\circ
\partial_i$.
\end{definition}
The operator $\theta_i$ can be thought of as a type of fission, where a
graph splits up into
$i$ particles. See Figure~\ref{ngon2}.

Extend $\theta_i$ to 
$$\theta_i: S(P\G^{\mathcal O})\to S(P\G^{\mathcal O})$$
as a derivation. Notice that $\theta_1=\partial_G=\phi_1$.

\begin{theorem}\label{whoknows}
The following identities hold:
\begin{itemize}
\item[a)] $\forall i  \neq j\hspace{5pt}
\theta_i\theta_j+\theta_j\theta_j=0$
\item[b)] $\forall i \hspace{5pt}
\theta_i^2=0$
\end{itemize}
\end{theorem}
\emph{[Proof]}

We prove a). Statement b) is similar. We show that
$$\theta_i\theta_j +\theta_j\theta_i \colon \G^{\mathcal O}\to
S^{i+j}(\G^{\mathcal O})$$ is zero, which is enough. If the $i$-gon and
$j$-gon attach to two different sets of edges, they can be applied in
either order to get the same (unoriented) result. Keeping track of
orientation, one sees that they anticommute.  Attaching one disk, and
then the other to an edge of the original disk is the same as adding a
bigger disk with an ordered pair of two sides marked for collapse.  We
may now apply our analysis from the proof of Lemma 1 to show that the
terms cancel in pairs.
$\hfill
\Box$

\begin{corollary}
Let $I$ be a subset of $\mathbb N$, finite or infinite. Let
$\theta_I = \sum_{i\in I}\theta_i$. Then $\theta_I^2=0$.
\end{corollary}

\begin{proposition}
$\theta_i$ is canonically zero at the level of homology.
\end{proposition}

\emph{[Proof]}

That $\theta_i$ is compatible with homology follows since
$\theta_1=\partial_G$.

A similar argument to Proposition~\ref{ich} shows that $\theta_i$ vanishes
canonically on homology.
$\hfill\Box$
\vskip5pt

The operator $\theta_i$ can be defined for disconnected graphs as well,
as we alluded to earlier. Suppose we start with a graph with $k$ connected
components. A $2i$-gon attaches to one of these and it fissions into
$i$ components. In order to get a well defined map, the remaining
$k-1$ components must be distributed with the
$i$ fission components in all possible ways, which leads to more
complicated formulas.

\section{Compatibility}

It is unclear if there is a theory of $Lie_\infty$ bialgebras; a search of
MathSciNet yields no hits. 
Under some obvious generalizations of the definition of Lie bialgebra to
the case of higher order operations on the symmetric algebra,
 the higher degree fusion operations are not
compatible with the higher degree fission operations.
Interestingly, degree 2 fission is compatible with degree 2 fusion
 on the
subcomplex of connected graphs with no separating edges. As was noted in
\cite{cv} this is not the case on the full complex $\mathcal G^{\mathcal
O}$.
\begin{definition}
Let $P^{irred}\mathcal G^{\mathcal O}$ be the subcomplex of $\mathcal
G^{\mathcal O}$ spanned by connected (primitive) graphs with no
separating edges (irreducible).
\end{definition}

\begin{theorem}\label{compat} On $P^{irred}\mathcal G^{\mathcal O}$ the
following equation holds:
$$\theta_2\phi_2(X\odot Y)+\phi_2(\theta_2(X)\odot Y)+
(-1)^{x}\phi_2(X\odot \theta_2(Y))=0.$$

\end{theorem}

\emph{[Proof]}

The bracket $\phi_2$ and cobracket $\theta_2$ coincide with the operations
$[\cdot,\cdot]$ and $\theta$ defined in \cite{cv} for the commutative
operad. In that paper, we defined everything in terms of contracting pairs
of half-edges, but the operations are easily seen to match. (In fact, we
mentioned a ``dotted line notation" in that paper which is very close to
the definition of
$\phi_2$ considered here.)
Now use the argument from \cite{cv} Theorem 1, which holds even if the
vertices are labelled by the operad $\mathcal O$.
$\hfill \Box$

\end{document}